\documentclass[graybox]{svmult}

\usepackage{type1cm}        
\usepackage{makeidx}         
\usepackage{graphicx}        
\usepackage{multicol}        
\usepackage[bottom]{footmisc}

\usepackage{newtxtext}       %
\usepackage[varvw]{newtxmath}       
\usepackage{epstopdf}
\usepackage{caption}
\usepackage{subcaption}

\newcommand{\ha}{\frac{1}{2}}

\newcommand{\tF}{\tilde{F}}

\def\be{\begin{equation}}
\def\ee{\end{equation}}
\newcommand{\ei}[0]{\end{itemize}}
\newcommand{\beann}[0]{\begin{eqnarray*}}
\newcommand{\eeann}[0]{\end{eqnarray*}}
\def\bea{\begin{eqnarray}}
\def\eea{\end{eqnarray}}
\def\ba{\begin{array}{l}\displaystyle}
\def\ea{\end{array}}

\makeindex             

\begin{document}

\title{Boundary effects on wave trains in the Exner model of sedimentat transport}

\author{E. Macca and G. Russo}
\institute{Emanuele Macca \at Universit{y} of Catania, \email{emanuele.macca@phd.unict.it}
\and Giovanni Russo \at Universit{y} di Catania, \email{giovanni.russo1@unict.it}
}

%
\maketitle

\abstract{The aim of this work is to study the effect of imposing different types of boundary conditions in the Exner model describing shallow water equations with sedimentation, when a train of waves is imposed at one of the boundaries \cite{TesiPhD,MaccaExner}. 
The numerical solver is a second order finite-volume scheme, with semi-implicit time discretization based on Implicit-Explicit (IMEX) schemes, which guarantee 
 better stability properties than explicit ones. We show the effect of spurious reflected waves at the right edge of the computational domain,  propose two remedies, and show how such spurious effects can be reduced by suitable non-reflecting boundary conditions.}

\section{Introduction} 
Sediment transport in channels and rivers can be effectively described in the framework of the Saint-Venant model of shallow water equations \cite{SaintVenant_1,SaintVenant_2}, by the so called Exner model
\cite{Fernandez,CastroNieto,Hudson,LiuBeljadid,MURILLO20108704,RzadkiewiczMariotti,VanRijn}, which, in one space dimension, constitutes a system of three partial differential equations for the water depth $h$, the discharge $q$ and the sediment thickness $z_b$, the first two equations describing conservation of water height and flux, while the last equation, describing the sediment evolution, is obtained by closing the system with a suitable constitutive relation which links the sediment flux to the water depth averaged velocity $u=q/h$.

The governing equations have the structure of a $3\times 3$ quasiliner strictly hyperbolic system, whose eigenvalues $\lambda_1<\lambda_2<\lambda_3$ represent the propagation velocities of the three families of waves. 

The flows we are interested in are {\em subcritical\/}, i.e.\ the water speed is always smaller than the propagation speed of the surface waves.  

Finite volume schemes have been successfully adopted  in computing the numerical solution of such system \cite{Castro-FernandezNieto2,CastroNieto,BonaventuraMultilayer}. 

For weak interaction, which depends on the coupling between the sediment and the water velocity, the wave speed $|\lambda_2|$ corresponding to the sediment transport is much smaller that the speeds of the surface waves, $|\lambda_1|$ and $|\lambda_2|$.

Under such conditions, the system can be considered {\em stiff\/}.
Explicit schemes have to satisfy the classical CFL condition on the time step, based on the fastest eigenvalue of the system. 

For such a reason, semi-implicit schemes have been proposed, which treat implicitly the fast water waves, thus requiring a much less restrictive CFL condition. 

Recently, Garres-Díaz et al. (2022) proposed a semi-implicit $\theta-$method approach for sediment transport models \cite{GarresFernandez} by which, choosing $\theta>\ha$ in the semi-implicit method, an increase in both efficiency and stability is obtained \cite{CasulliCattani}.

In \cite{MaccaExner} the authors show that in many cases the evolution of the sediment can be captured without resolving the detailed evolution of the surface water waves. 

Because of the slow evolution of the sediment, one has to integrate the equations for times which are usually much longer than the travel time of surface waves in the computational domain. It would be therefore desirable to impose non-reflecting boundary conditions on such waves, in order to avoid spurious reflections at the boundary. 

The goal of the paper is to study the effect of different boundary conditions on the right edge of the computational domain for  long time simulation, when a wave train is imposed on the left edge of the computational domain. 




The paper is structured in the following way: in Section 2 the one-dimensional model equations are introduced and discussed; in Section 3 the details of the semi-implicit numerical method are presented. In Section 4 different numerical right boundary treatment are shown. In Section 5 we present several tests to numerically assess the effectiveness of the various boundary conditions in presence of right traveling wave trains on free-surface and some numerical strategies to reduce the reflecting surface waves. Finally, in Section 6 we draw some conclusion.

\section{1D Exner Model}
	Let us start from 1D shallow water equations on a time independent bathymetry $b(x)$:
	\begin{equation}
    \label{Shallow_water}    
    \begin{cases}
        h_t + q_x =0\\
     \displaystyle   q_t + \left(\frac{q^2}{h} + \frac{g}{2}h^2 \right)_x = -ghb_x,
    \end{cases}
    \end{equation}
    \begin{figure}[!ht]
	\hspace{-1.cm}	
	\includegraphics[scale = 0.39]{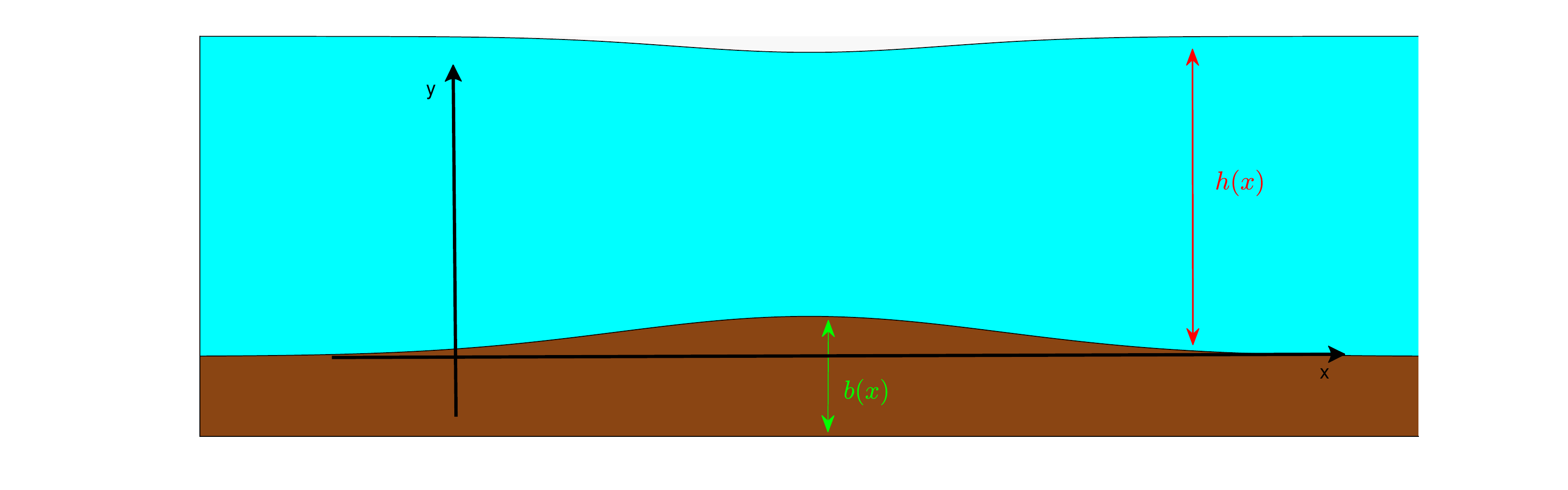}
	\vspace{-0.4 cm}
	\caption{Shallow water equations: water-flow $h(x)$ and bottom topography $b(x).$}
	\label{Fig_Ex_1}
\end{figure}
where $x$ denotes the space coordinate along the axis of the channel and $t$ is time; $q(x,t)$  represents the water flux per unit width (discharge) and $h(x,t)$ the water thickness; $g$ the acceleration due to gravity;  $b(x)$ denotes the bottom topography; furthermore, the following relation holds $q(x,t) =h(x,t)u(x,t),$ in which $u$ is the depth average horizontal velocity as displayed in Figure \ref{Fig_Ex_1}.

Replacing $q$ by $hu$ in the second equation, and making use of the equation for $h$, and assuming always positive $h$, one obtains the following equation for the velocity:
\begin{equation}
\label{u_eq}
 u_t + u u_x + g h_x = 0
\end{equation}

\begin{figure}[!ht]
	\hspace{-1cm}
	\includegraphics[scale = 0.39]{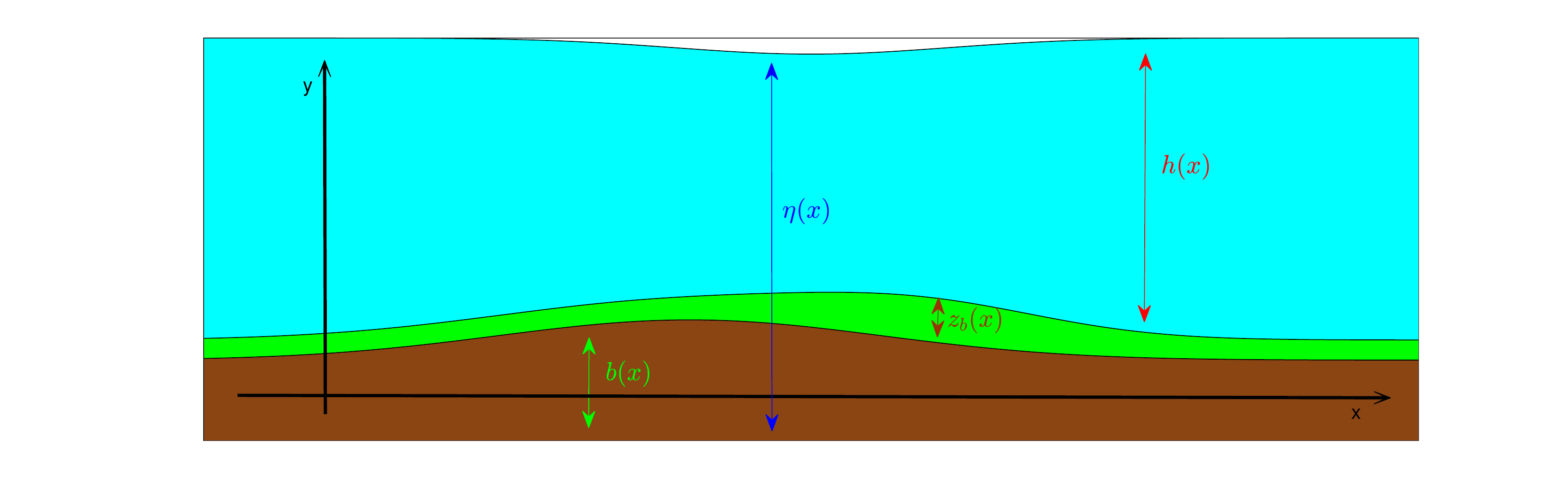}
	\vspace{-0.4 cm}
	\caption{1D Exner model: water surface $\eta(x);$ water-flow $h(x);$ sedimental layer $z_b(x)$ and bottom topography $b(x).$}
	\label{Fig_Ex_2}
\end{figure}

The Exner system, adopted in this work, has been obtained coupling shallow water equation \eqref{Shallow_water} with an equation describing the sediment transport. 
The system then can be written as 
\begin{equation}
	    \label{Ex_sis}
	    \begin{cases}
	        h_t + q_x = 0,\\
	        \displaystyle q_t + \left(\frac{q^2}{h} + \ha gh^2\right)_x = -gh(b + z_b)_x, \\
	        (z_b)_t + (q_b)_x = 0.  
        \end{cases}
\end{equation}
where the driving force on the RHS now depends on the total bottom height (bathymetry $b(x)$ and sediment layer thickness $z_b(x,t)$, see Figure \ref{Fig_Ex_2}). 

In order to close system \eqref{Ex_sis} we adopt the Grass formulation  \cite{Grass,QianLi,CastroNieto} 
\begin{equation}
	    \label{q_b}
	    q_b = \xi A_g u|u|^{m-1}
\end{equation}
in which $m$ comes out from experimental results, $m\in[1,4],$ $A_g$ is the interaction between the fluid and the sediment, in general $A_g\in]0,a[$ with $a>0$ and $\xi = {1}/{(1-\rho_0)}$ where $\rho_0$ is the porosity of the sediment layer. Throughout this paper we shall assume that the porosity is constant.

System \eqref{Ex_sis} does not admit a conservation form, because the bottom profile $S(x,t) = b(x) + z_b(x,t)$ contains the unknown field. 
In terms of the unknown field 
\[
	    W = ( h, q, S)^\top
\]
the system can be written in quasilinear form as


\begin{equation}
	    \label{non_cons_term_1}
	    \frac{\partial W}{\partial t} + A(W)\frac{\partial W}{\partial x} = 0,
\end{equation}
where,
\[
    A(W) = \begin{bmatrix} 0 & 1 & 0 \\gh - u^2 & 2u & gh \\ \alpha & \beta & 0     \end{bmatrix}, \> {\rm with}
\]
and $
\alpha \equiv \partial_h q_b, \> \beta \equiv \partial_q q_b$. 
Assuming $u>0$ in the whole domain one has $\beta = m \xi A_g u^{m-1}/h$ and $\alpha = -u\beta$.

The non-conservative system is strictly hyperbolic if and only if the characteristic polynomial: 
\[
    p_{\lambda}(\lambda) = -\lambda((u-\lambda)^2 -gh) + gh\beta(\lambda -u)
\]
has three distinct real roots $\lambda_1<\lambda_2<\lambda_3.$

As $A_g \rightarrow0,$ $\beta$ vanishes, and the three eigenvalues become $\lambda_1 = u-c$, 
$\lambda_2 = 0$ and $\lambda_3 = u+c$, with 
$c=\sqrt{gh}$. 

The main interest of this paper is the treatment of boundary conditions for the Exner model adopting a semi-implicit IMEX scheme. The IMEX methods present advantages in both stability and efficiency in regimes for which the local Froude number $F_r = |u|/c$ is relatively small, say $F_r<1/2$, and the quantity $\beta$ is much smaller than 1. For this reason, we only focus on small value of the interaction coefficient  $A_g$ between the fluid and the sediment.

For sufficiently small values of $A_g$, such that $\beta\ll 1$, performing an asymptotic expansion of the eigenvalues, we obtain (assuming $u>0$) 
\begin{align}
\label{lamd_1}
    \lambda_1 &= u - \sqrt{gh} - \beta\frac{\sqrt{gh}}{2(1 - F_r)} + O(\beta^2) \\ \label{lamd_2}
    \lambda_2 &= \frac{\beta g h u}{gh - u^2} + O(\beta^2) = \beta u /(1-F_r^2) + O(\beta^2) \\ \label{lamd_3}
    \lambda_3 &= u + \sqrt{gh} + \beta\frac{\sqrt{gh}}{2(1 + F_r)} + O(\beta^2)   
\end{align}

In the next section we shall derive the semi-implicit scheme that we adopt in the paper, in which, under the assumption of small Froude number, we treat implicitly the surface waves and explicitly the sediment wave.

To obtain automatically a well-balanced method, it is convenient to  rewrite the Exner system \eqref{Ex_sis} in terms of  $\eta(x,t) = h(x,t) + b(x) + z_b(x,t) $, which  represents the elevation of the undisturbed water surface, in place of the water thickness $h$ (see Figure \ref{Fig_Ex_2}). System \eqref{Ex_sis} becomes:
	\begin{equation}
	    \label{Ex_sis_eta}
	    \begin{cases}
	        \eta_t + (q+q_b)_x = 0\\
	        q_t + (qu)_x + gh(\eta)_x = 0 \\
	        (z_b)_t + (q_b)_x = 0  
        \end{cases}
	\end{equation}
 with $h = \eta-b-z_b$, $u=q/h$, and $q_b$ given by Eq. \eqref{q_b}.
 
\section{Semi implicit scheme} \label{1D semi implicit}
	The introduction and the derivation of the semi-implicit scheme will be presented in this section. Since we are interested in the preservation of the sediment evolution we adopt an implicit treatment of the surface water waves, while the slow wave corresponding to the sediment evolution is treated explicitly. In particular, we focus on semi-implicit schemes, first and second order accurate, in space and time.

    \subsection{Mesh}
In this work we consider uniform meshes in 1D. Space domain is represented by a segment $\Omega$ and it is uniformly discretized in $N$ cells, each of size $\Delta x$.  

The center of cell $\omega_i=[x_{i-1/2};x_{i+1/2}]$ is denoted by $x_i$, $i=1, \ldots, N$.

The time domain $\mathcal{T}=[0,T]$ with final time $T>0$ is split into time intervals $[t^n,t^{n+1}]$, $n\in \mathbb{N}$, and time-steps $\Delta t_n=t^{n+1}-t^n$ subject to a CFL (Courant-Friedrichs-Lewy) like condition\footnote{Although it is better to assign $\Delta t$ dynamically at each time step by imposing some CFL condition, here we shall adopt a constant time step $\Delta t$, in order to simplify the notation in the description of the method.}.
	
	Definitely, we denote by $U_i^n$ an approximation on the mean value of $U$ over cell $I_i$ at time $t=t^n,$
	$$U^n_i\cong \frac{1}{\Delta x}\int_{x_{i-\ha}}^{x_{i+\ha}}U(x,t^n)dx. $$
	
	\subsection{First order scheme}
	Let us consider system \eqref{Ex_sis_eta}. Following the idea proposed in \cite{Boscarino-Filbet}, system \eqref{Ex_sis_eta} can be rewritten as a large system of ODE's, in which we adopt suitable discrete operators for the approximation of space derivatives. The key point in \cite{Boscarino-Filbet} is to identify which specific term has to be treated implicitly and which can be treated explicitly. 
    
    Following \cite{Boscarino-Filbet}, we rewrite system \eqref{Ex_sis_eta} identifying the terms that will be treated explicitly, and the ones that require an implicit treatment: 
    \begin{equation}
        \label{ode_bosca_semi}
        U' = \tilde{K}(U_E,U_I)
    \end{equation}
    where $U = [\eta, q, z_b]^\top$ and $\tilde{K}(U_E,U_I)$ is given by
    \begin{equation}
        \label{ode_form_semi}
        \tilde{K}(U_E,U_I) = 
        \begin{bmatrix}
            -(q_I+(q_b)_E)_x \\
            -((qu)_E)_x - gh_E(\eta_I)_x \\
            -((q_b)_E)_x
        \end{bmatrix}
    \end{equation}
in which, the subscript $E$ and $I$ denote which term has to be treated explicitly and which implicitly. 

By replacing differential operators with suitable discrete ones, with a slight abuse of notation the semi-discrete scheme can be written in the form
    \begin{equation}
        \label{ode_bosca}
        U' = K(U_E,U_I).
    \end{equation}
with
    \begin{equation}
        \label{ode_form_part}
        K(U_E,U_I) = 
        \begin{bmatrix}
           -\hat{D}_x((q_B)_E)  & -D_x(q_I) \\ -\hat{D}_x((qu)_E) &  -gh_ED_x(\eta_I) \\ -\hat{D}((q_b)_E) &
        \end{bmatrix}.
    \end{equation}
	Finally, the fully-discrete first order in time semi-implicit scheme can be written as:
	\begin{equation}
	    \begin{cases}
	        \label{First order}
	        \eta^{n+1} = \eta^n - \Delta t \hat{D}_x(q_b^n) - \Delta t D_x(q^{n+1}),\\
	        q^{n+1} = q^n - \Delta t \hat{D}_x(q^nu^n) - \Delta tgh^nD_x(\eta^{n+1}), \\
	        z_b^{n+1} = z_b^n - \Delta t \hat{D}_x(q_b^n),
	    \end{cases}
	\end{equation}
	where the discrete operators $D_x$ and $\hat{D}_x$ applied to a given flux function $F(U)$ are respectively defined as: 
	\begin{itemize}
	    \item $D_x(\tF_i) = \frac{\tF_{i+\ha} - \tF_{i-\ha}}{\Delta x},$ in which $\tF_{i\pm\ha}$ is suitably defined on cell edges (see Remark \ref{rem1} below); 
	    \item $\hat{D}_x(\tF_i) = \frac{\tF_{i+\ha} - \tF_{i-\ha}}{\Delta x},$ where  $\tF_{i+\ha} = \ha\Bigl( \tF(U_{i+\ha}^{-}) + \tF(U_{i+\ha}^{+}) - \alpha_{i+\ha}\big(U_{i+\ha}^{+} - U_{i+\ha}^{-}\big)\Bigr)$ is the Rusanov flux and $\alpha_{i+\ha}$ is related to the eigenvalues of the explicit sub system. In our case  $\alpha_{i+\frac12} = \max(|u^-_{i+\ha}|,|u^+_{i-\ha}$ and is much smaller than $|\lambda_1|$ and $|\lambda_3|$. $U_{i\pm\ha}^{\pm}$ are obtained by piecewise linear reconstruction with  MinMod slope limiter.
	\end{itemize}
 In order to distinguish explicit part from implicit one in system \eqref{First order}, we set $\eta^*$ and $q^*$ the explicit part of first and second equation. In this way, the system can be rewritten as:
	\begin{equation}
	    \begin{cases}
	        \label{First order star}
	         q^{*} = q^n - \Delta t \hat{D}_x(q^nu^n);\\
	         \eta^{*} = \eta^n - \Delta t \hat{D}_x(q_b^n) - \Delta t \hat{D}_x(q^{*});\\
	        \eta^{n+1} = \eta^* + g\Delta t^2 D_x(h^nD_x(\eta^{n+1}));\\
	        q^{n+1} = q^* - \Delta tgh^nD_x(\eta^{n+1}); \\
	        z_b^{n+1} = z_b^n - \Delta t \hat{D}_x(q_b^n); \\ 
            h^{n+1} = \eta^{n+1} - z_b^{n+1} - b.
	    \end{cases}
	\end{equation}

 \begin{remark}
 \label{rem1}
     The differential operators $D_x$, in the third and fourth lines of equation \eqref{First order star} respectively, are so discretized: 
     $$
        \left.D_x(h^nD_x(\eta^{n+1}))\right|_i =  \frac{1}{\Delta x^2}\Bigl(h_{i+\ha}^n(\eta_{i+1}^{n+1} - \eta_{i}^{n+1}) - h_{i-\ha}^n(\eta_{i}^{n+1} - \eta_{i-1}^{n+1})\Bigr),  $$
         where $h_{i\pm\ha}^n = \ha(h_{i}^n + h^n_{i\pm1});$ and
         $$ \left. D_x(\eta^{n+1})\right|_i = \frac{1}{2\Delta x}( \eta_{i+1}^{n+1} - \eta_{i-1}^{n+1}).$$
 \end{remark}

\subsection{Second order scheme}
    In the previous section we derived a semi-implicit scheme which is first order in time and second order in space. In this section we focus on the second order IMEX Runge-Kutta procedure \cite{AVGERINOS2019278,PareschiRusso} to obtain second order both in space and time. 

    In order to obtain a second order semi-implicit time discretization we adopt the technique introduced in \cite{Boscarino-Filbet}, which is based on (apparently) doubling system \eqref{ode_bosca}, and applying to this double system an IMEX Runge-Kutta for partitioned systems. 

An IMEX scheme is defined by a double Butcher tableau of the form 
\begin{equation*}
\label{IMEX_tableau}
    \begin{array}{c|cc}
            \tilde{c} & \tilde{A} \\ \hline
                      & \tilde{b}^\top
        \end{array} 
        \hspace{3 cm}
    \begin{array}{c|cc}
            {c} & {A} \\ \hline
                      & {b}^\top
        \end{array} 
\end{equation*}
where the lower triangular matrix $\tilde{A}\in\mathbb{R}^{s\times s}$ with zero diagonal and the vectors $\tilde{c}$ and $\tilde{b}$ characterize the explicit part of the scheme, while 
the triangular matrix ${A}\in\mathbb{R}^{s\times s}$ and the vectors $c$ and $b$ identify the implicit part of the IMEX scheme.  
    If the implicit scheme is stiffly-accurate, i.e.\ if $a_{si} = b_i,\> i=1,\ldots,s$ then the last stage of the method coincides with the numerical solution. 
    Applying the IMEX scheme to system \eqref{ode_bosca_semi} one obtains


    \begin{itemize}
        \item Stage values: For $i=1,\ldots,s$ compute
        \begin{align*}
            U^{(i)}_E & = U^n + \Delta t\sum_{j=1}^{i-1}a_{i,j}^EK\left(U_E^{(j)},U_I^{(j)}\right)\\
            U^{(i)}_I & = U^n + \Delta t
            \left(\sum_{j=1}^{i-1}a_{i,j}^I K\left(U_E^{(j)},U_I^{(j)}\right) + a_{i,i}^I K\left(U_E^{(i)},U_I^{(i)}\right)\right).
        \end{align*}
        \item Numerical solution:\\
        $U^{n+1} = U_I^{(s)}.$
    \end{itemize}
    For more details about the method the reader may consult \cite{Boscarino-Filbet}.


Here we consider the IMEX scheme defined by the following  double Butcher tableau \cite{Boscarino-Filbet}: 
\begin{equation}
\label{tableau}
    \begin{array}{c|cc}
             & 0 &  \\
            c & c & 0\\ \hline
            & 1-\gamma & \gamma
        \end{array} 
        \hspace{3 cm}
        \begin{array}{c|cc}
            \gamma & \gamma &  \\
            1 & 1-\gamma & \gamma\\ \hline
            & 1-\gamma & \gamma
        \end{array}
\end{equation}
with $\gamma = 1 - \frac{1}{\sqrt{2}}$ and $c = \frac{1}{2\gamma}.$ Applying the scheme defined by \eqref{tableau} we have:
    \begin{enumerate}
        \item $U_E^{(1)} = U^n;$
        \item $U_I^{(1)} = U^n + \Delta t\gamma K(U_E^{(1)},U_I^{(1)});$
        \item $U_E^{(2)} =  U^n + \Delta tc K(U_E^{(1)},U_I^{(1)});$
        \item $U_I^{(2)} = U^n + \Delta t(1-\gamma) K(U_E^{(1)},U_I^{(1)}) + \Delta t\gamma K(U_E^{(2)},U_I^{(2)});$
        \item $U^{n+1} = U_I^{(2)}.$
    \end{enumerate}
    \begin{remark}
Let observe that $U_E^{(2)},U_I^{(2)}$ and $U_I^{(1)}$ have a common term, thus step 3 and 4 may be rewritten as:
    \begin{align*}
        U_E^{(2)} &= (1-\frac{c}{\gamma})U^n + \frac{c}{\gamma}U_I^{(1)};  \\
        U_I^{(2)} &= (1-\frac{1-\gamma}{\gamma})U^n + \frac{1-\gamma}{\gamma}U_I^{(1)} +\Delta t\gamma K(U_E^{(2)},U_I^{(2)}).
    \end{align*} 
    \end{remark}

\section{Right boundary treatment}
As we have mentioned before, imposing a wave train on the left boundary of the domain, when the signal arrives to the right edge,  the right boundary condition is crucial in determining the wave pattern in the domain, after a long time. 

Here we compare three different boundary conditions:
\begin{enumerate}
    \item Simple transmissive conditions, i.e.\ zero Neumann conditions on all field variables (NC). 
    \item Simple wave approximation (SC). The solution on the last cell is assumed to be approximated by a flat bottom and a simple wave for shallow water equations corresponding to the largest eigenvalue. 
    \item Absorbing layer conditions (AC). A layer is introduced to the right of the computational domain, in which a damping is imposed on the equations, which damps the waves once they pass the right edge.


    In the latter approach, we split the entire domain $\Omega$ into two regions $\Omega = \Omega_E \cup \Omega_A,$ where $\Omega_E$ is the domain in which we adopt the undamped Exner model, while $\Omega_A$ represents the absorbing domain that we adopt to damp the waves. In practice, for the first order scheme we apply Eq. \eqref{First order star} in $\Omega_E$, while in the domain $\Omega_A$ we adopt the following scheme
    \begin{equation}
	    \begin{cases}
	        \label{First order star mod}
	         q^{*} = q^n - \Delta t \hat{D}_x(q^nu^n) - \Delta t(q^n-q^0)\varphi(x);\\
	         \eta^{*} = \eta^n - \Delta t \hat{D}_x(q_b^n) - \Delta t \hat{D}_x(q^{*}) - \Delta t(\eta^n-\eta^0)\varphi(x);\\
	        \eta^{n+1} = \eta^* + g\Delta t^2 D_x(h^nD_x(\eta^{n+1}));\\
	        q^{n+1} = q^* - \Delta tgh^nD_x(\eta^{n+1}); \\
	        z_b^{n+1} = z_b^n - \Delta t \hat{D}_x(q_b^n) - \Delta t(z_b^n-z_b^0)\varphi(x),
	    \end{cases}
	\end{equation}
 where $\varphi(x) := \Bigl((x-x_R)/\sigma(l)\Bigr)^2,$ with $x_R$ we denote the right boundary of the domain $\Omega_E$ and $\sigma(l)$ represents a length scale which is chosen to be much larger than the wavelength of the wave train in order to avoid reflected waves.
\end{enumerate}

\section{Numerical experiments}
This section focuses on the numerical comparison between the 3 different techniques for treating the condition at the right edge. We consider a single numerical experiment in which we show snapshots of of the numerical solution at various times, obtained with the three techniques to delineate the different effects.

In all our tests we adopt the following value for the time step:
\[
    \Delta t_n = {\rm CFL}\frac{\Delta x}{\lambda_{\rm max}^n }
\]
where $\lambda_{\rm max}^n$ is the maximum, over cells, of the spectral radius of absolute of matrix $A(U^n)$. We adopt an approximation of the eigenvalues of $A$ given in \eqref{lamd_1}-\eqref{lamd_2}-\eqref{lamd_3}, neglecting  the $O(\beta^2)$ terms.  
Furthermore we monitor the material CFL
$$
    \textrm{MCFL} = \frac{u^n_{\rm max} \Delta t_n}{\Delta x} 
$$
where $u^n_{\rm max} = \max_j{|u^n_j|}$, and verify that it is always less than 1.

\begin{figure}[!ht]
    \vspace{-3cm}
    \centering
	\includegraphics[width=0.7 \textwidth]{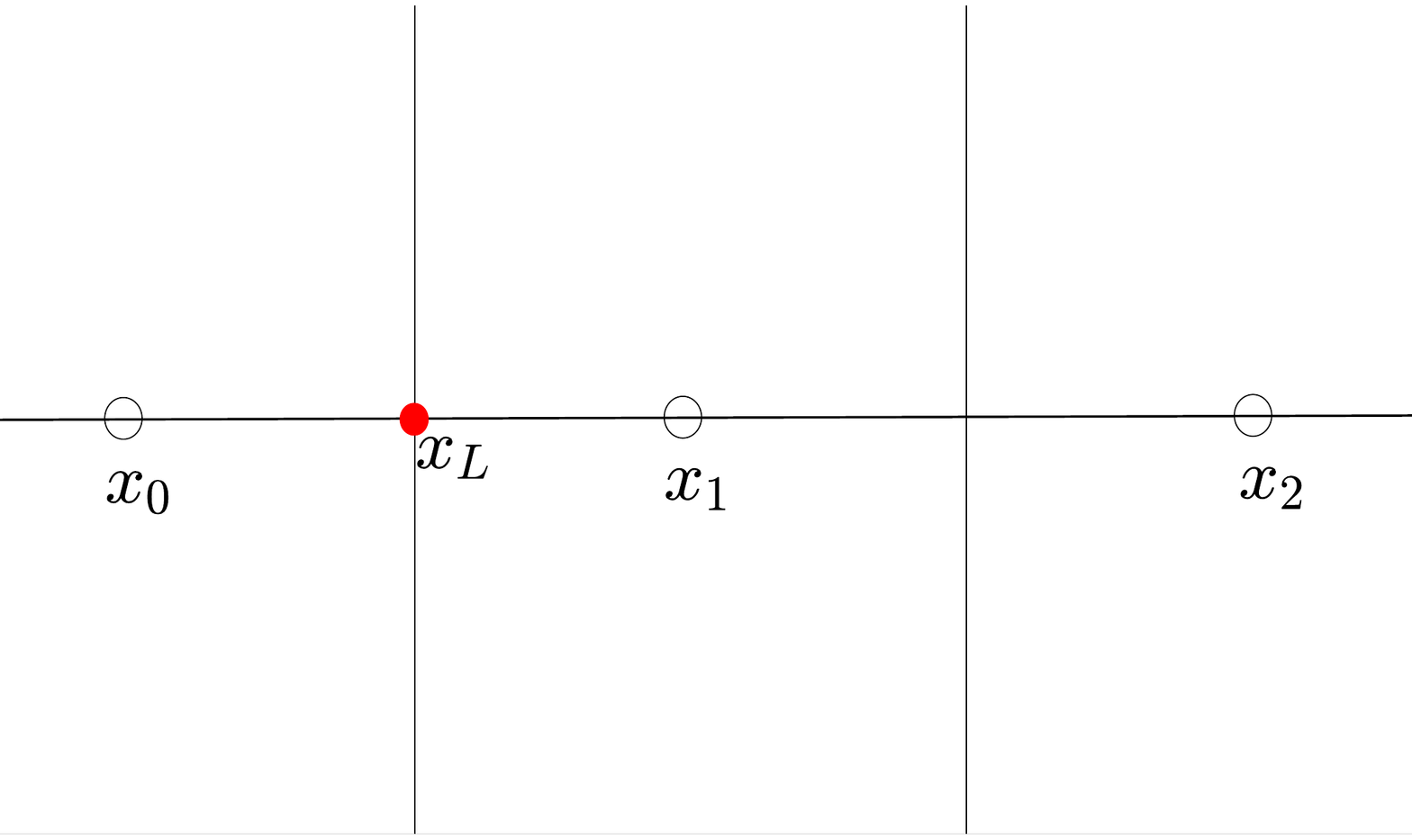}
	\vspace{-3cm}
	\caption{Left-boundary mesh with ghost point 0.}
	\label{Fig_bound}
\end{figure}

\begin{figure}[!ht]
	\hspace{-1cm}
	\includegraphics[width=1.15 \textwidth]{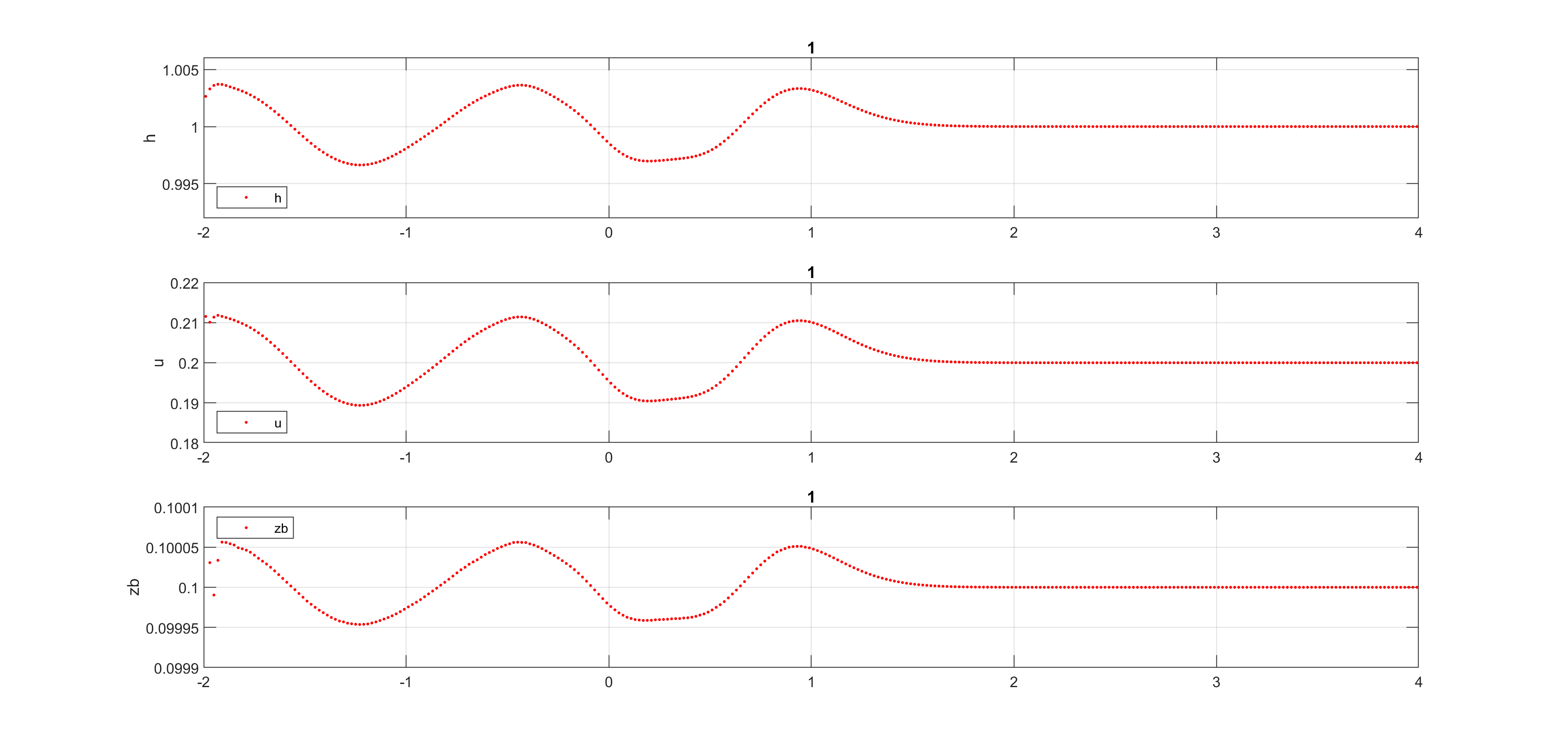}
    \vspace{-0.6cm}
	\caption{1D Exner model: water-flow $h;$ velocity $u$; and sedimental layer $z_b$ where the Neumann zero condition (NC) has been adopted for the right boundary treatment at final time $t = t_{\rm fin}^0 = 1.$}
	\label{Fig_free_0}
\end{figure}
\begin{figure}[!ht]
	\hspace{-1cm}
	\includegraphics[width=1.15 \textwidth]{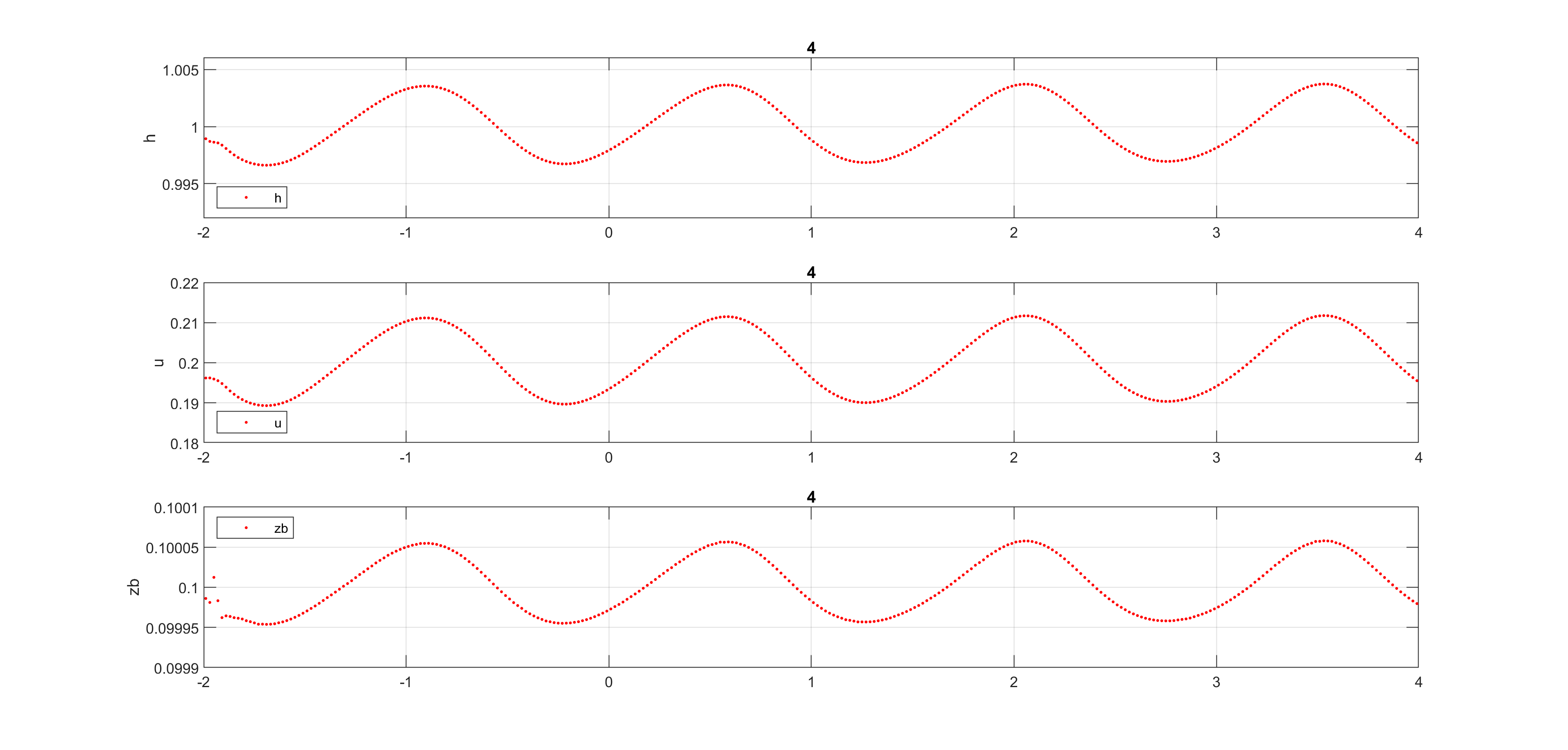}
	\vspace{-0.6cm}
	\caption{1D Exner model: water-flow $h;$ velocity $u$; and sedimental layer $z_b$ where the Neumann zero condition (NC) has been adopted for the right boundary treatment at final time $t = t_{\rm fin}^1 = 4.$}
	\label{Fig_free_1}
\end{figure}
\begin{figure}[!ht]
	\hspace{-1cm}
	\includegraphics[width=1.15 \textwidth]{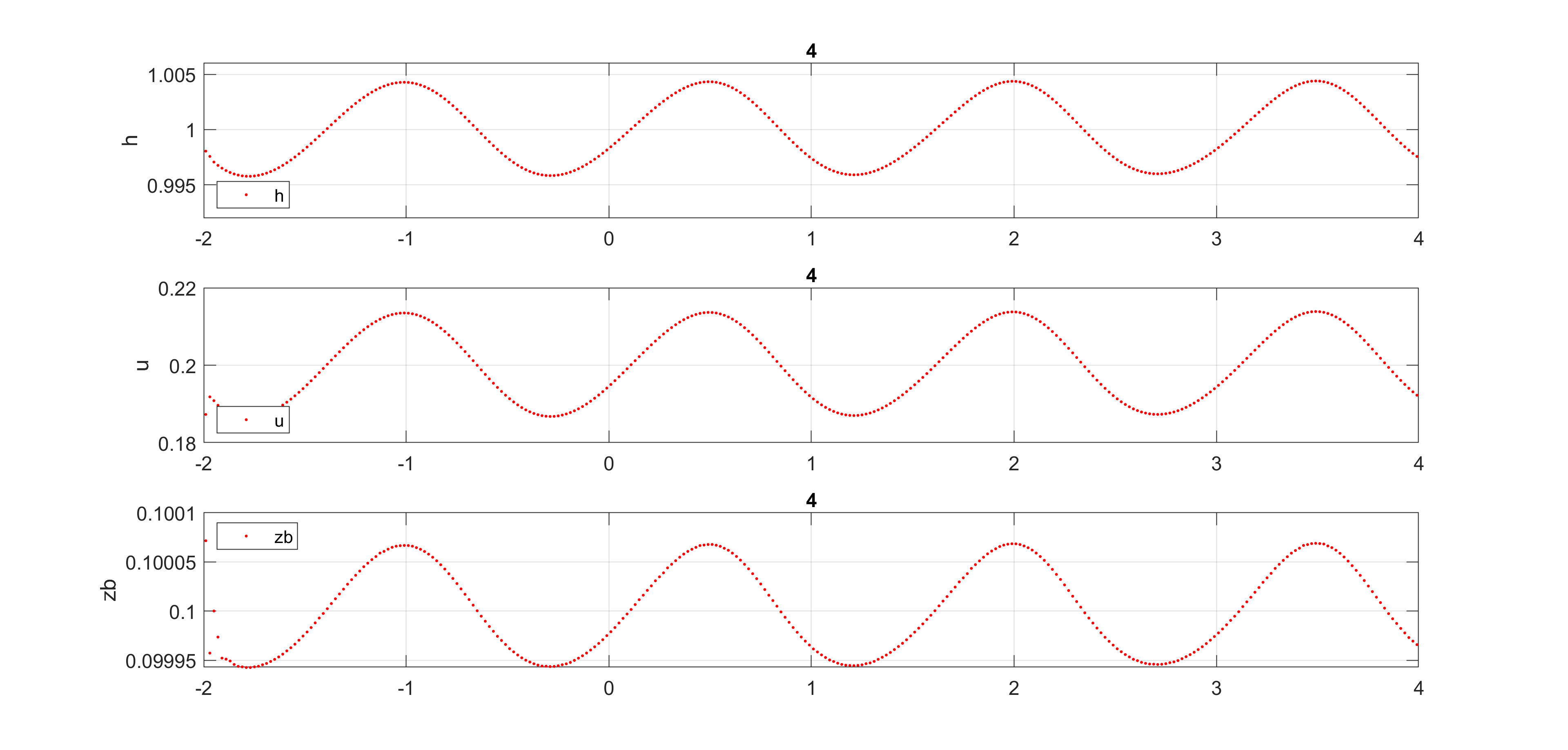}
	\vspace{-0.6cm}
	\caption{1D Exner model: water-flow $h;$ velocity $u$; and sedimental layer $z_b$ where the Burgers approximation (SC) has been adopted for the right boundary treatment at final time $t = t_{\rm fin}^1 = 4.$}
	\label{Fig_bu_1}
\end{figure}
\begin{figure}[!ht]
	\hspace{-1cm}
	\includegraphics[width=1.15 \textwidth]{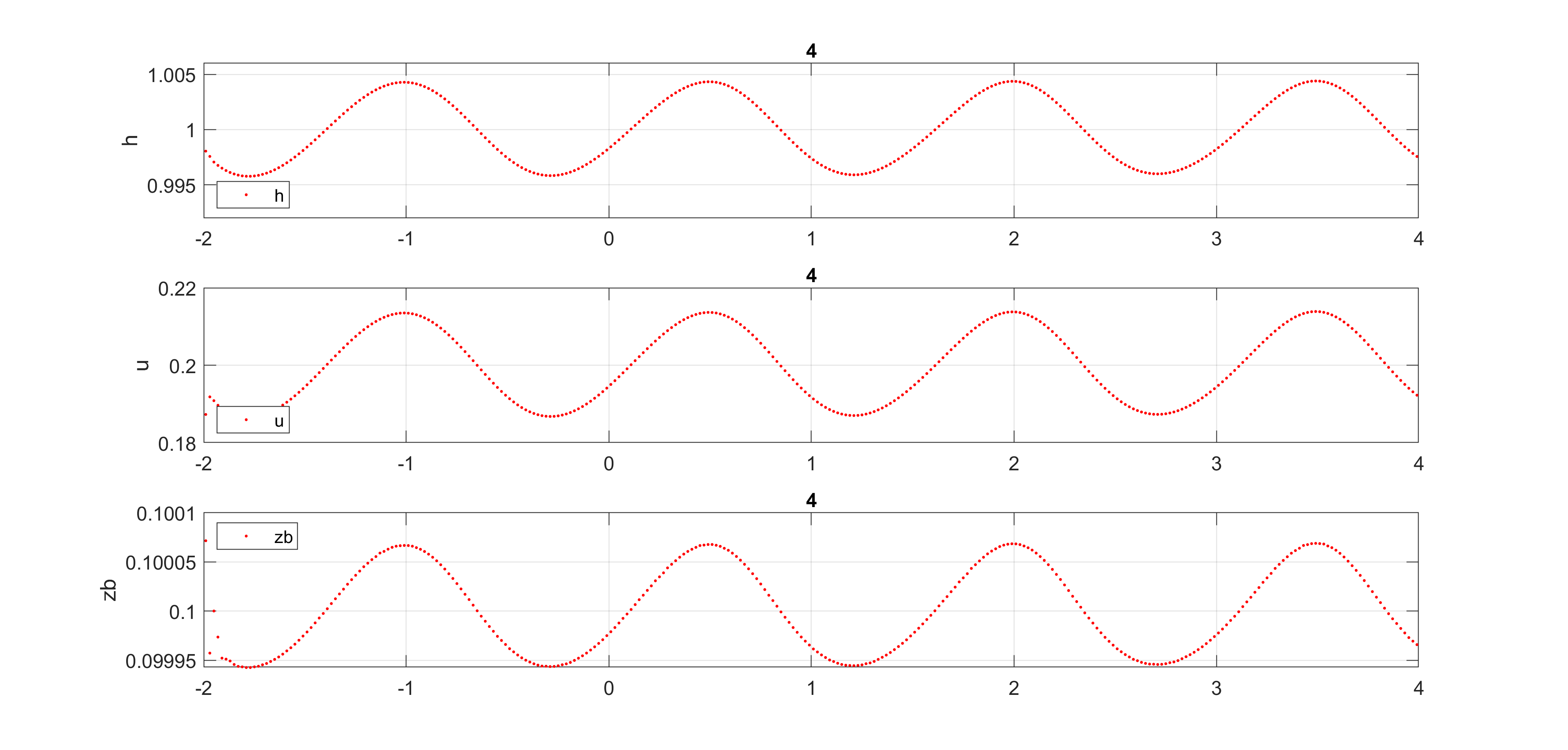}
	\vspace{-0.6cm}
	\caption{1D Exner model: water-flow $h;$ velocity $u$; and sedimental layer $z_b$ where the absorbing domain $\Omega_A$ (AC) has been adopted for the right boundary treatment at final time $t = t_{\rm fin}^1 = 4.$}
	\label{Fig_sm_1}
\end{figure}
Let us consider the one-dimensional Exner system \eqref{Ex_sis_eta} with the Grass equation \eqref{q_b} and the second order semi-implicit method illustrated before.
The common parameters are so set:  computational the domain is $\Omega = \Omega_E \cup \Omega_A = [-2,4]\cup[4,10];$ CFL is set to $7.7,$ (which gives  $\textrm{MCFL} < 0.5$ in All tests); 
the Grass parameter IS $A_g = 0.1,$ $\xi = 1/(1-\rho_0)$ where $\rho_0 = 0.2$ and $m=3;$ the acceleration due to the gravity IS $g = 9.81;$ the bottom topography $b$ is set to zero.
We show snapshots of the solution at four different times: 
a) before the waves touch the right boundary of the computational domain $\Omega_E$,   $t_{\rm fin}^0 = 1$;  after a reflected wave reached the left boundary, $t_{\rm fin}^1 = 4$, after two full reflections, $t_{\rm fin}^2 = 8$ and after a long time, $t_{\rm fin}^3 = 10^6.$
\begin{figure}[!ht]
	\hspace{-1cm}
	\includegraphics[width=1.15 \textwidth]{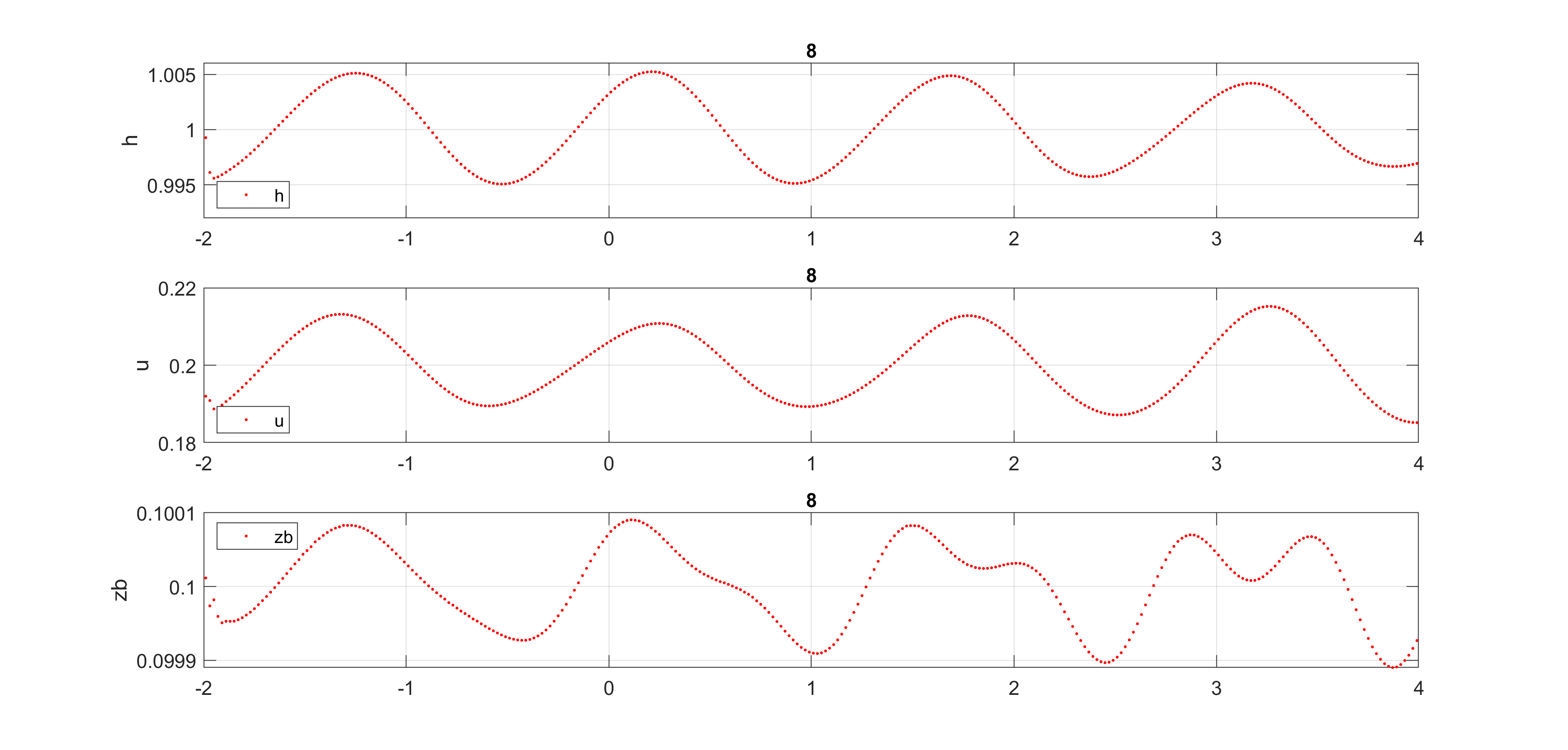}
	\vspace{-0.6cm}
	\caption{1D Exner model: water-flow $h;$ velocity $u$; and sedimental layer $z_b$ where the Neumann zero condition (NC) has been adopted for the right boundary treatment at final time $t = t_{\rm fin}^2 = 8.$}
	\label{Fig_fr_2}
\end{figure}
\begin{figure}[!ht]
	\hspace{-1cm}
	\includegraphics[width=1.15 \textwidth]{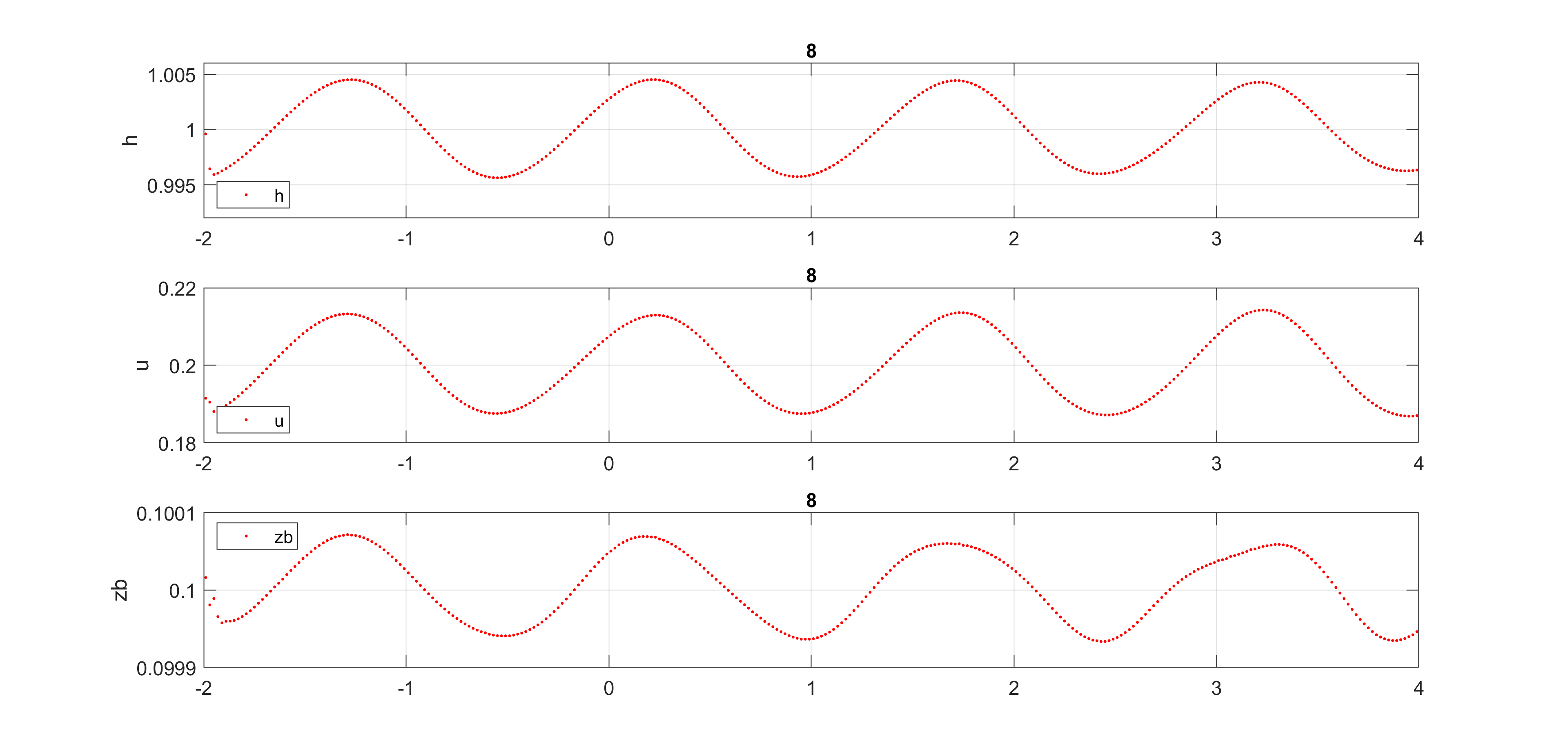}
	\vspace{-0.6cm}
	\caption{1D Exner model: water-flow $h;$ velocity $u$; and sedimental layer $z_b$ where the simple wave approximation (SC) has been adopted for the right boundary treatment at final time $t = t_{\rm fin}^2 = 8.$}
	\label{Fig_bu_2}
\end{figure}
\begin{figure}[!ht]
	\hspace{-1cm}
	\includegraphics[width=1.15 \textwidth]{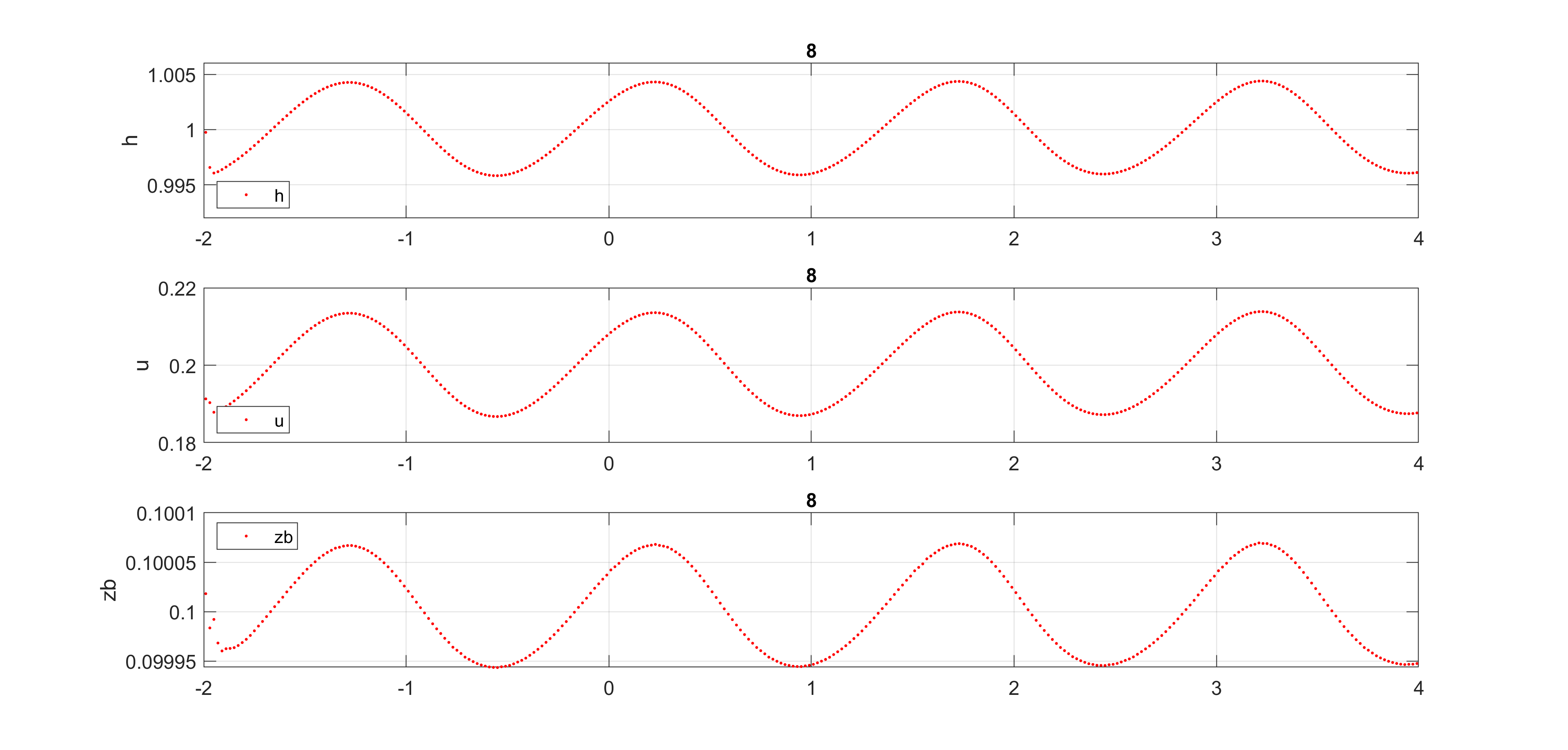}
	\vspace{-0.6cm}
	\caption{1D Exner model: water-flow $h;$ velocity $u$; and sedimental layer $z_b$ where the absorbing domain $\Omega_A$ (AC) has been adopted for the right boundary treatment at final time $t = t_{\rm fin}^2 = 8.$}
	\label{Fig_sm_2}
\end{figure}
 The initial conditions are 
\[
    (h(x,0),u(x,0),z_{b}(x,0)) = (1,0.2,0.1).
\]

 The boundary conditions on the left edge are imposed by assigning the velocity at the entrance, $u(x_L,t) = \phi(t)$, and imposing a compatibility condition on $h$ 
 given by Eq.~\eqref{u_eq}, assuming that 
 the sediment profile is  be flat at the entrance, i.e. $\partial_x z_b(0,t) = 0$. 
 All these conditions are imposed to second order accuracy on the discrete scheme by assigning the following value at the ghost cell immediately to the left of the domain $\Omega_E$ (see Fig.\ \ref{Fig_bound}).

 $$
\begin{bmatrix}
      h_0 \\ u_0 
    \end{bmatrix} = \begin{bmatrix}
       h_1 + 1/g(\phi_t\Delta x + 0.5((u_1)^2 - \phi^2)) \\ 2\phi-u_1 
\end{bmatrix} 
$$

 In all tests we adopt  $\phi(t) = 0.2 + \mathcal{A}(\sin(\omega t))$ in which $\mathcal{A}$ and $\omega$ denote amplitude and frequency of the waves in our case set to 0.01 and 14 respectively \cite{TesiPhD}, and $\phi_t \equiv d\phi/dt.$

\begin{figure}[!ht]
	\hspace{-1cm}
	\includegraphics[width=1.15 \textwidth]{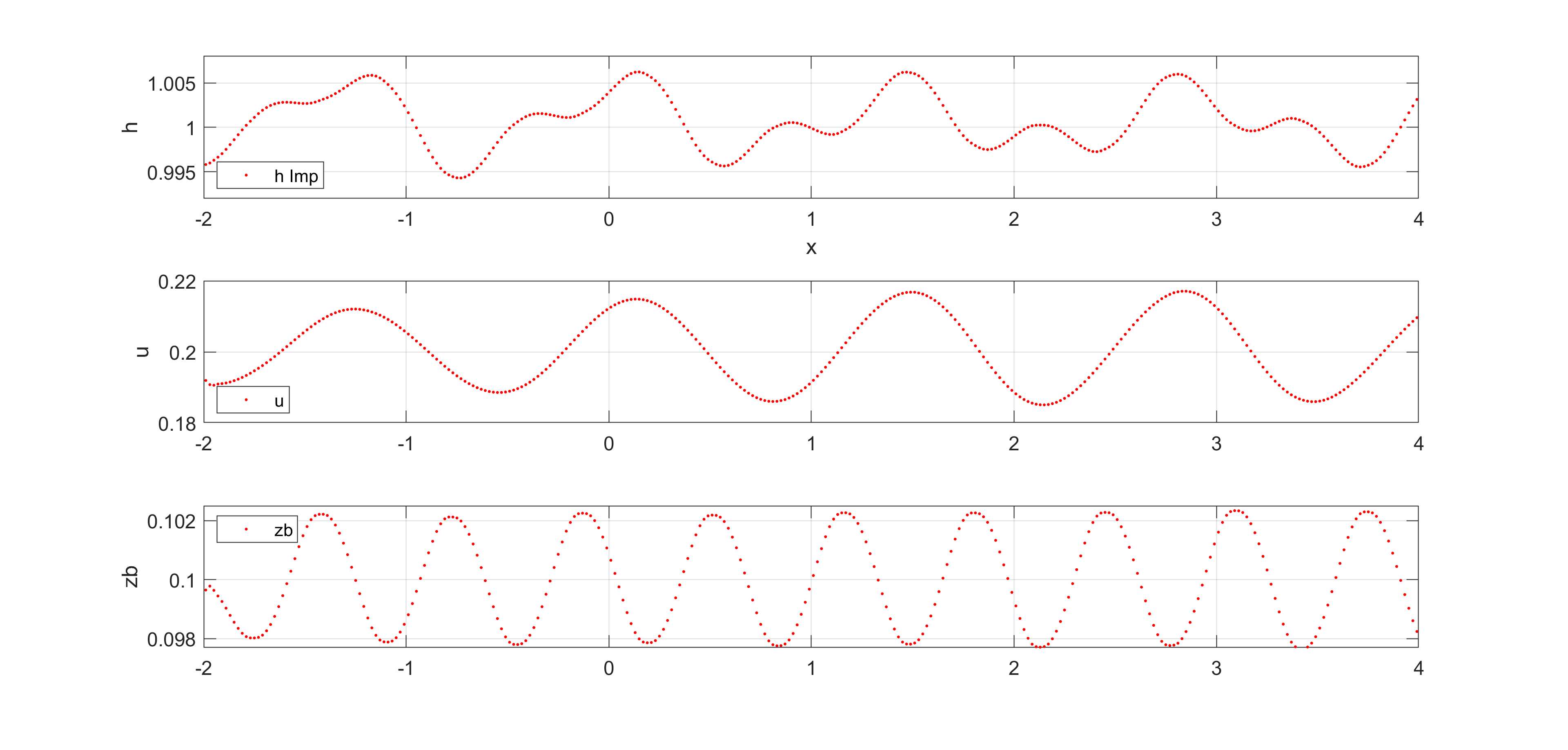}
	\vspace{-0.6cm}
	\caption{1D Exner model: water-flow $h;$ velocity $u$; and sedimental layer $z_b$ where the Neumann zero condition (NC) has been adopted for the right boundary treatment at final time $t = t_{\rm fin}^3 = 10^6.$}
	\label{Fig_fr_3}
\end{figure}
\begin{figure}[!ht]
	\hspace{-1cm}
	\includegraphics[width=1.15 \textwidth]{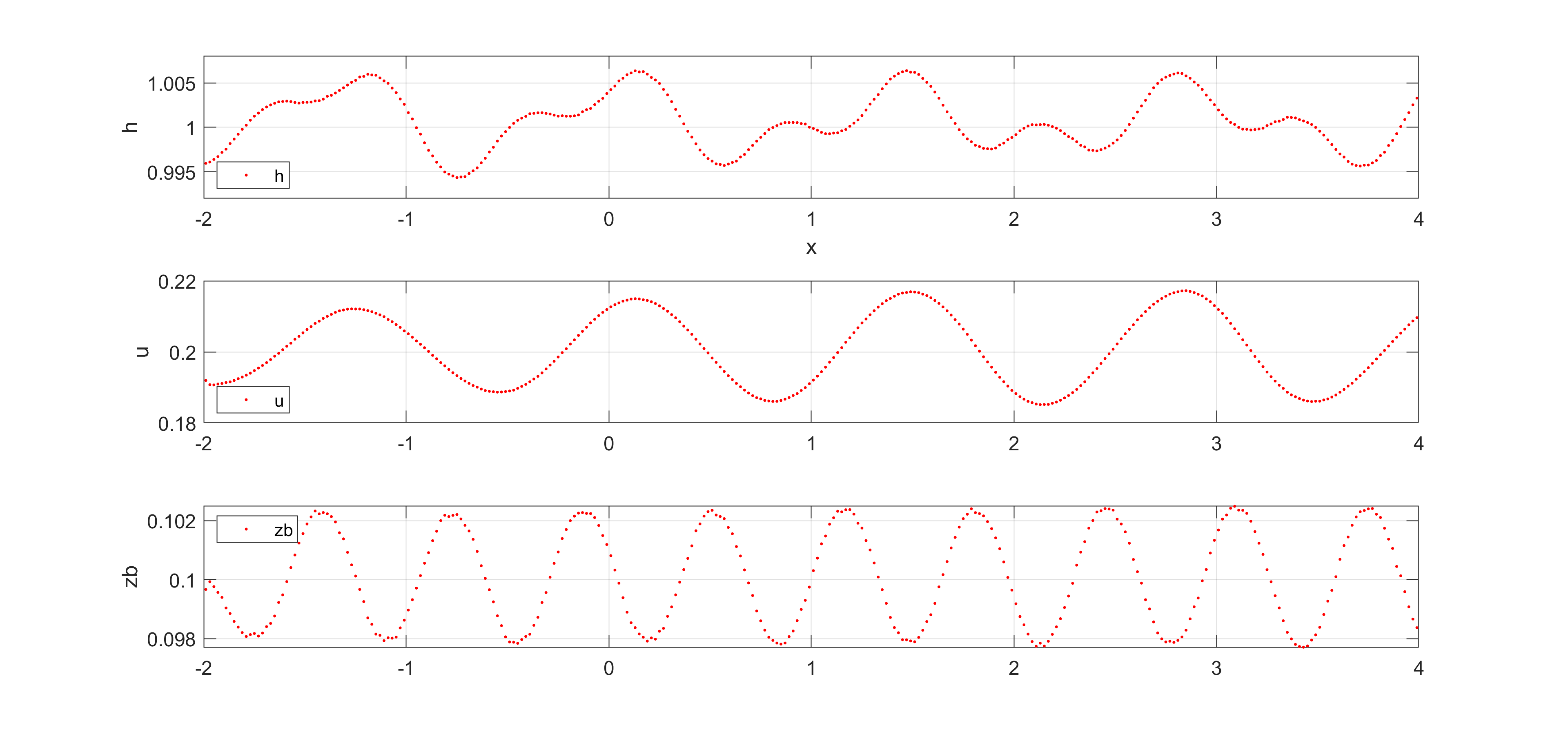}
	\vspace{-0.6cm}
	\caption{1D Exner model: water-flow $h;$ velocity $u$; and sedimental layer $z_b$ where the simple wave approximation (SC) has been adopted for the right boundary treatment at final time $t = t_{\rm fin}^3 = 10^6.$}
	\label{Fig_bu_3}
\end{figure}
\begin{figure}[!ht]
	\hspace{-1cm}
	\includegraphics[width=1.15 \textwidth]{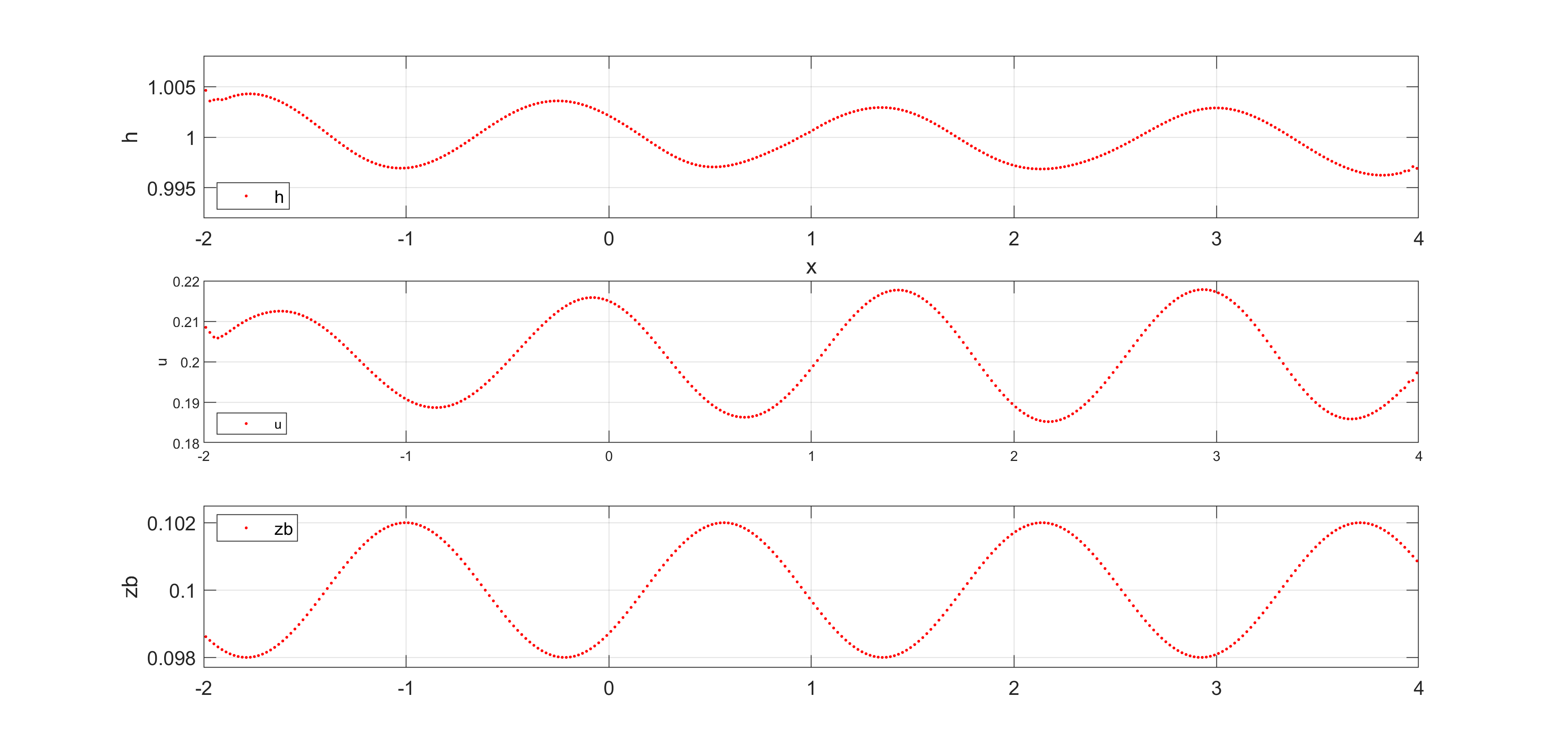}
	\vspace{-0.6cm}
	\caption{1D Exner model: water-flow $h;$ velocity $u$; and sedimental layer $z_b$ where the absorbing domain $\Omega_A$ (AC) has been adopted for the right boundary treatment at final time $t = t_{\rm fin}^3 = 10^6.$}
	\label{Fig_sm_3}
\end{figure}

 All the figures show the numerical solutions for water-flow $h;$ velocity $u$; and sedimental layer $z_b$ in which different conditions have been applied to the right boundary of the computational domain $\Omega_E.$ In particular, Figure \eqref{Fig_free_0} shows the numerical solutions before the waves touch the right boundary of $\Omega_E$; Figures \eqref{Fig_free_1}-\eqref{Fig_sm_1} exhibit the numerical solutions at short time $t=4$; Figures \eqref{Fig_fr_2}-\eqref{Fig_sm_2} display the numerical solutions at medium time $t=8;$ and Figures \eqref{Fig_fr_3}-\eqref{Fig_sm_3} show the numerical solutions at long time $t=10^6.$ Is is remarkable, from Figures \eqref{Fig_fr_3}-\eqref{Fig_sm_3} how, for a long time simulation, the effects of the reflecting waves produced by strategy (NC) and (SC), even if reduced, produce almost the same pattern (i.e. a period doubling in the sediment profile). 
 Analyzing Figures \eqref{Fig_free_1}-\eqref{Fig_sm_1} one observes that the effect of the small reflected wave is barely noticeable in all cases. 
 After two full reflections (Figures \eqref{Fig_fr_2}-\eqref{Fig_sm_2}), the use of the approximate simple wave near the right boundary mitigates the effect of the reflected wave. Finally, after a very long time, only the technique based on the absorbing layer is able to effective eliminate the reflected waves.
 It is peculiar to observe that the two approaches (NC) and (SC) after a long time give very similar results, i.e.\  the sediment profile assumes a periodic structure with a period which is twice the original one, while the $h$ profile  almost identical in the two cases, in spite of the fact that the effect of the reflected waves is considerably smaller for (SC) boundary conditions than in the case of (NC) ones.

\section{Conclusion}
The purpose of this work the comparison of three different techniques for the treatment of the right boundary condition when a wave train  is imposed on the left boundary for the one-dimensional Exner model \cite{TesiPhD,MaccaExner}.

Sediment transport occurs at a much lower speed than surface waves. Consequently, in order to follow the motion of the sediment over a large fraction of the computational domain, the equations have to be integrated times which are much longer than surface waves travel time. If no particular care is taken at the boundaries, during such a a long time span, surface waves reflect back and forth inside the domain, thus completely changing the sediment dynamics that one would observe by integrating the same equations on a much longer computational. domain, i.e. with no reflections.

In order to reduce computational time, it is of paramount relevance to adopt non-reflecting boundary conditions on the right side of the domain. 

Three different numerical boundary conditions have been compared. 
 The first is the classic zero Neumann condition on all the unknowns (NC). This technique, extremely simple in the implementation, produced considerable reflective effects, which completely compromise the solution after few reflections. 
 
 The second technique is to use, in the ghost cell immediately at the right of the computational domain $\Omega_E$, a numerical solution approximated by the a simple wave of shallow water equations with constant bathymetry (SC). Such a technique has been introduced in the PhD thesis \cite{TesiDuven} in the context of Euler equations of compressible gas dynamics, where it was shown to be effective reducing the impact of reflected waves by approximately one order of magnitude. Nevertheless, for a large number of reflections, this technique is not sufficient to make the reflection effects negligible.
 
 The third technique consists in adding a sufficiently large auxiliary domain,  in which the waves are gradually dampened to the point of disappearing without creating reflexive phenomena (AC). This technique completely solves the problem introduced by the first two techniques paying a higher computational cost due to the extension of the computational domain.

There are still a few things that require improvements and generalizations:
\begin{itemize}
\item Two approximations, adopted in the present paper, may reduce the effectiveness of the technique based on the approximate simple waves. Such approximations consists in assuming that the inflow sediment profile is flat on the left edge of the domain, so that we can easily impose the compatibility condition on the water depth based on Eq.\ \eqref{u_eq}, and in assuming that we can approximate a simple wave of the whole Exner system by a simple wave for shallow water on constant bathymetry. Both approximations have been made in order to simplify the treatment and resort so analytical expressions. It would be interesting to explore  semi-analytical conditions which may improve the boundary treatment by (SC). 
\item It would be interesting to further understand the period-doubling effect on the sediment profile induced by the reflected waves, which appear to be very robust and almost independent on the detail of the reflection mechanism. 
\item In the current method the stability condition on the time step is determined by the fluid velocity, which may be much larger than the sediment wave speed. We plan to construct a scheme in which the CFL condition that determines the stability is based on the sediment wave rather than on the fluid velocity.
\item To explore other, more efficient techniques to eliminate reflected waves. 
\end{itemize}
All such generalizations and extension will be subject of future investigation.

\section*{Acknowledgements} 
 This research has received funding from the European Union’s NextGenerationUE – Project: Centro Nazionale HPC, Big Data e Quantum
Computing, “Spoke 1” (No. CUP E63C22001000006). E. Macca was partially supported by GNCS No. CUP E53C22001930001 Research Project “Metodi numerici per problemi differenziali multiscala: schemi di alto ordine, ottimizzazione, controllo”. E. Macca and G.Russo would like to thank the Italian Ministry of Instruction, University and Research (MIUR) to support this research with funds coming from PRIN Project 2017 (No. 2017KKJP4X entitled “Innovative numerical methods for evolutionary partial differential equations and applications”). E. Macca and G. Russo are members of the INdAM Research group GNCS.

\bibliographystyle{plain}
\bibliography{biblio}

\end{document}